\documentclass[11pt,a4paper]{amsart}
\usepackage{graphicx,multirow,array,amsmath,amssymb}

\newtheorem{theorem}{Theorem}

\begin{document}

\title{Equilateral stick number of knots}
\author[H. Kim]{Hyoungjun Kim}
\address{Department of Mathematics, Korea University, Anam-dong, Sungbuk-ku, Seoul 136-701, Korea}
\email{Kimhjun@korea.ac.kr}
\author[S. No]{Sungjong No}
\address{Department of Mathematics, Korea University, Anam-dong, Sungbuk-ku, Seoul 136-701, Korea}
\email{blueface@korea.ac.kr}
\author[S. Oh]{Seungsang Oh}
\address{Department of Mathematics, Korea University, Anam-dong, Sungbuk-ku, Seoul 136-701, Korea}
\email{seungsang@korea.ac.kr}

\thanks{This research was supported by Basic Science Research Program through
the National Research Foundation of Korea(NRF) funded by the Ministry of Science,
ICT \& Future Planning(MSIP) (No.~2011-0021795).}

\begin{abstract}
An equilateral stick number $s_{=}(K)$ of a knot $K$ is defined to be the minimal number
of sticks required to construct a polygonal knot of $K$ which consists of equal length sticks.
Rawdon and Scharein \cite{RS} found upper bounds for the equilateral stick numbers of
all prime knots through 10 crossings by using algorithms in the software KnotPlot.
In this paper, we find an upper bound on the equilateral stick number of a nontrivial knot $K$
in terms of the minimal crossing number $c(K)$ which is $s_{=}(K) \leq 2c(K) + 2$.
Moreover if $K$ is a non-alternating prime knot, then $s_{=}(K) \leq 2c(K) - 2$.
Furthermore we find another upper bound on the equilateral stick number for composite knots
which is $s_{=}(K_1 \sharp K_2) \leq 2c(K_1) + 2c(K_2)$.
\end{abstract}

\maketitle

\section{Introduction} \label{sec:intro}
A knot is a closed curve in $3$-space $\mathbb{R}^3$.
Knots appear in many physical systems in the natural sciences,
for example, in molecular chains such as DNA and proteins,
and they have been considered to be useful models for structural analysis of these molecules.
A knot can be embedded in many different ways in $3$-space, smooth or piecewise linear.
By a {\em polygonal knot\/} we will mean a knot which consists of finite line segments,
called {\em sticks\/}, glued end-to-end.
This representation of knots is very useful for many applications
because microscopic level molecules are more similar to rigid sticks than flexible ropes.

Define the {\em stick number\/} $s(K)$ of a topological knot type $K$ as
the minimal number of sticks required to realize $K$ as a polygonal knot.
Randell \cite{R} first showed that a nontrivial polygonal knot needs at least six sticks and
$s(3_1)=6$ and $s(4_1)=7$.
In addition, all five and six crossing prime knots $5_1$, $5_2$, $6_1$, $6_2$ and $6_3$,
the square and granny knots and two eight crossing knots $8_{19}$ and $8_{20}$ are
only known knots having stick number 8.
Several upper and lower bounds on the stick number for various classes of knots were founded.
Negami \cite{N} used graph theory to provide upper and lower bounds for the stick number
in terms of the crossing number $c(K)$, $(5+ \sqrt{8c(K)+9})/2 \leq s(K) \leq 2c(K)$.
Later Huh and Oh \cite{HO} improved Negami's upper bound to $s(K) \leq 3(c(K)+1)/2$.
Furthermore $s(K) \leq 3c(K)/2$ for a non-alternating prime knot $K$.
Also Calvo \cite{Ca} improved this lower bound to $(7 + \sqrt{8c(K)+1})/2 \leq s(K)$.
For the stick number of a composite knot, Adams, Brennan, Greilsheimer and Woo \cite{ABGW}
and Jin \cite{J} independently proved that $s(K_1 \sharp K_2) \leq s(K_1) + s(K_2) -3$.

\begin{figure}[h]
\begin{center}
\includegraphics{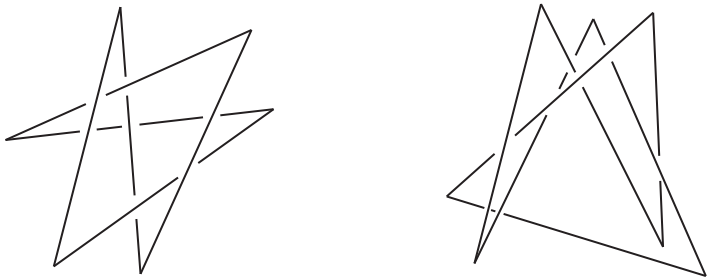}
\end{center}
\caption{Equilateral knots of $3_1$ and $4_1$}
\label{fig1}
\end{figure}

In this paper we are interested in another quantity concerning polygonal knots.
A {\em equilateral knot\/} is a polygonal knot which consists of equal length sticks.
The {\em equilateral stick number\/} $s_{=}(K)$ of a knot $K$ is defined to be
the minimal number of equal length sticks required to construct
a equilateral knot representation of $K$.
Little is known about the equilateral stick number.
Adams et al. \cite{ABGW} proved that $s_{=}(K_1 \sharp K_2) \leq s_{=}(K_1) + s_{=}(K_2)$.
Also Rawdon and Scharein \cite{RS} used algorithms in the software KnotPlot to compute
upper bounds for the equilateral stick numbers of all prime knots through 10 crossings.
They showed that all such knots except seven knots can be constructed by the same
number of equal length sticks as their stick numbers.
Figure \ref{fig1} shows projections of the trefoil and the figure-eight knots
constructed by minimal number of equal length sticks.

In this paper we establish an upper bound on the equilateral stick number of nontrivial knots
in terms of the minimal crossing number.
Furthermore we find another upper bound on the equilateral stick number of composite knots.

\begin{theorem} \label{thm1}
Let $K$ be a nontrivial knot.
Then $s_{=}(K) \leq 2c(K) + 2$.
Moreover if $K$ is a non-alternating prime knot, then $s_{=}(K) \leq 2c(K) - 2$.
\end{theorem}

\begin{theorem} \label{thm2}
Let $K_1$ and $K_2$ be nontrivial knots.
Then $s_{=}(K_1 \sharp K_2) \leq 2c(K_1) + 2c(K_2)$.
Moreover if one of knot is non-alternating prime,
then $s_{=}(K_1 \sharp K_2) \leq 2c(K_1) + 2c(K_2) -4$,
and if both knots are non-alternating prime,
then $s_{=}(K_1 \sharp K_2) \leq 2c(K_1) + 2c(K_2) -8$.
\end{theorem}

Remark that the equilateral stick number of a knot with at most 10 crossings
is very close to its stick number as the result in \cite{RS}.
Thus the reader may curious about that we probably reduce the upper bound of the equilateral stick number
to $3c(K)/2 + n$ for some positive number $n$ by using the similar method in \cite{HO}.
But there is an obstacle to extending this method to constructing an equilateral realization.
In \cite{HO}, vertical line segments vary in length, so if we adjust them in a uniform length,
then unexpected crossing changes might happen.

\section{Arc index}

In this section we introduce another minimality invariant, an {\em arc index\/}.
There is an open-book decomposition of $\mathbb{R}^3$ which has open half-planes as pages
and the standard $z$-axis as the binding axis.
We may regard each page as a half-plane $H_{\theta}$ at angle $\theta$
when the $xy$-plane has a polar coordinate.
It can be easily shown that every knot $K$ can be embedded in an open-book decomposition
with finitely many pages so that it meets each page in a simple arc.
Such an embedding is called an {\em arc presentation\/} of $K$.
The {\em arc index\/} $a(K)$ is defined to be the minimal number of pages
among all possible arc presentations of $K$.
The left of Figure \ref{fig2} shows an arc presentation of the trefoil knot.
Here the points of $K$ on the binding axis are called binding indices,
assigned by $1,2, \cdots , a(K)$ from bottom to top.

We mention two rotating operations on arc presentations of knots.
Like turning pages of a book we can rotate the pages along the binding axis of an arc presentation.
This changing is called {\em page rotating\/}.
Similarly we can also rotate binding indices
so that for some fixed integer $n$ each binding index $i$ goes to $i+n$
where if $i+n > a(K)$, then we use $i+n - a(K)$ instead of $i+n$.
This is call {\em binding index rotating\/}.
Clearly both rotations do not change the knot type of $K$.

Bae and Park established an upper bound on the arc index in terms of the crossing number.
Corollary $4$ and Theorem $9$ in \cite{BP} provide that $a(K) \leq c(K)+2$,
and moreover $a(K) \leq c(K)+1$ if $K$ is a non-alternating prime knot.
Later Jin and Park improved the second part of Bae and Park's theorem.
Theorem $3.3$ in \cite{JP} provides that if $K$ is a non-alternating prime knot, then $a(K) \leq c(K)$.
Thus we have the following;

\begin{theorem} \label{thm:ac}
Let $K$ be any nontrivial knot. Then $a(K) \leq c(K)+2$.
Moreover if $K$ is a non-alternating prime knot, then $a(K) \leq c(K)$.
\end{theorem}

\section{Proof of Theorem \ref{thm1}}

\begin{figure}[h]
\begin{center}
\includegraphics{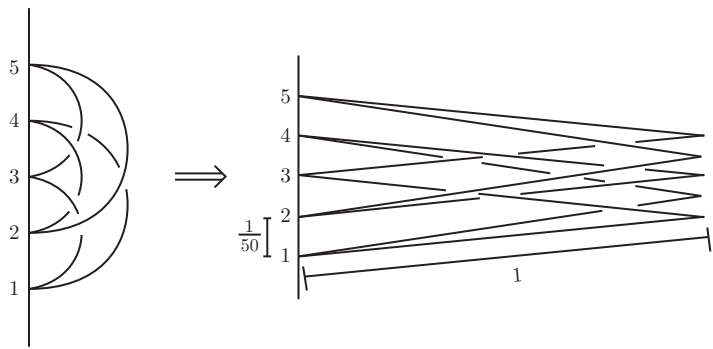}
\end{center}
\caption{Arc presentation and equilateral knot}
\label{fig2}
\end{figure}

In this section we prove Theorem \ref{thm1}.
Let $K$ be a nontrivial knot with $a(K) = n$ and $K^a$ be an arc presentation of $K$.
Now we modify $K^a$ so that consecutive binding indices are close enough
and each arc consists of two adjacent equilateral sticks on the related page.
We may assume that the distance between any two consecutive binding indices is $\frac{1}{10n}$
for convenience and all sticks have length 1 as the right in Figure \ref{fig2}.
Already we obtain a equilateral knot of $K$ consisting of $2 a(K)$ equal length sticks.

\begin{figure}[h]
\begin{center}
\includegraphics{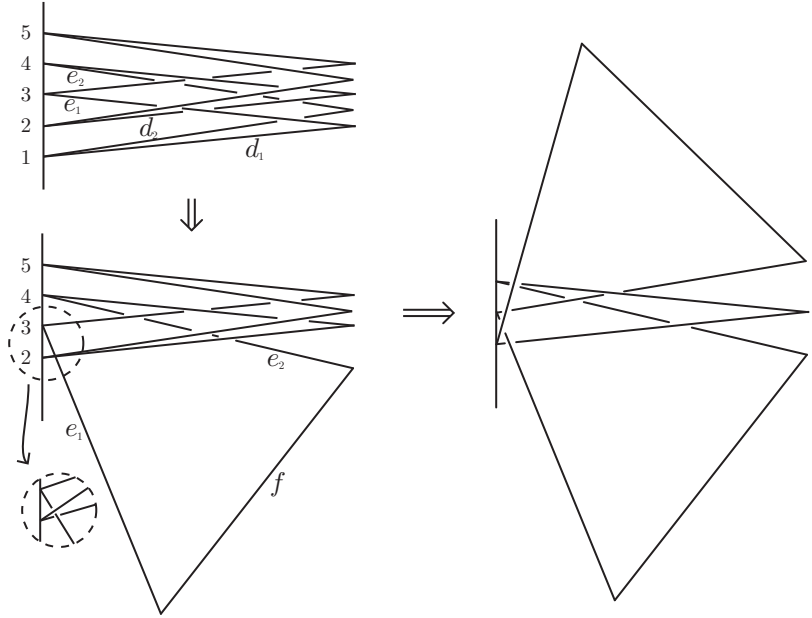}
\end{center}
\caption{Reducing some sticks}
\label{fig3}
\end{figure}

Now we will reduce two more sticks from $K^a$.
Let $d_1$ and $d_2$ be the two sticks adjacent to the binding index 1,
and $e_1$ and $e_2$ be the other two sticks adjacent to $d_1$ and $d_2$, respectively.
As in Figure \ref{fig3}, delete two sticks $d_1$ and $d_2$.
Rotate the stick $e_1$ clockwise along its endpoint at the binding axis
on the related page until $e_1$ is close enough to the binding axis.
Now rotate the other stick $e_2$ clockwise along its endpoint at the binding axis
on the related page so that the distance between the other endpoints of $e_1$ and $e_2$
which are not at the binding axis is 1.
Glue a new stick $f$ with length 1 to these endpoints.
Clearly this formation does not change the knot type,
but we reduce one stick from the original.

Note that there is no arc connecting the first and the last binding indices.
Otherwise, after some binding index rotating, this arc connects consecutive two binding indices.
This means that we can reduce the arc index, a contradiction.
Thus we can repeat the above argument via the two sticks $d'_1$ and $d'_2$ adjacent to the binding index $n$
and the other two sticks $e'_1$ and $e'_2$ adjacent to $d'_1$ and $d'_2$, respectively,
because these four sticks are different from the previous four sticks.
Thus we reduce one more stick again.

Hence the resulting equilateral knot consists of $2 a(K) -2$ sticks.
Theorem \ref{thm:ac} guarantees Theorem \ref{thm1}.

\section{Proof of Theorem \ref{thm2}}

\begin{figure}[h]
\begin{center}
\includegraphics{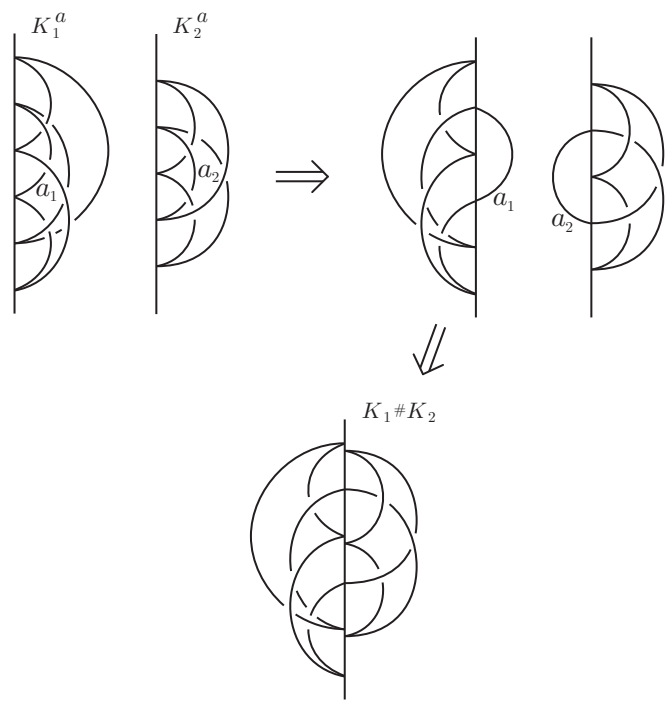}
\end{center}
\caption{Arc presentation of $K_1 \sharp K_2$}
\label{fig4}
\end{figure}

In this section we prove Theorem \ref{thm2}.
Let $K_1$ and $K_2$ be nontrivial knots
and $K^a_1$ and $K^a_2$ be their minimal arc presentations.
Since nontrivial knots have at least 5 arcs in their arc presentations,
we can choose arcs $a_1$ and $a_2$ of $K^a_1$ and $K^a_2$, respectively,
whose endpoints are not the first or the last binding indices of these arc presentations.
By using page rotating, we move all the arcs, except $a_1$, of $K^a_1$
and only $a_2$ of $K^a_2$ to the opposite sides as in Figure \ref{fig4}.
Then, delete two arcs $a_1$ and $a_2$,
and connect the two arcs of $K^a_1$ to the two arcs of $K^a_2$,
yielding a representation of $K_1 \sharp K_2$.
Note that the binding indices of $K^a_1$ except the connection points
must be apart from the binding indices of $K^a_2$.
Now we can apply the same arguments as in the proof of Theorem \ref{thm1}
independently to the arc presentations on both sides of the binding axis.
Thus we can reduce two sticks on both sides.

Hence the resulting equilateral knot consists of
$2(a(K_1)-1) + 2(a(K_2)-1) -4$ sticks.
Theorem \ref{thm:ac} guarantees Theorem \ref{thm2}.

\end{document}